\documentclass[12pt]{article}%
\usepackage{amsfonts}
\usepackage{amsmath,amssymb}
\usepackage[applemac]{inputenc}
\usepackage{amsmath,amssymb,fullpage}
\usepackage{color}
\usepackage{amsmath}
\usepackage{amssymb}
\usepackage{graphicx}%
\providecommand{\U}[1]{\protect\rule{.1in}{.1in}}

\newtheorem{definition}{Definition}[section]

\newtheorem{theorem}[definition]{Theorem}

\numberwithin{equation}{section}

\def\1B{\text{1\!\!I}}

\begin{document}

\date{20 March 2018}
\title{Reflected Advanced Backward Stochastic Differential Equations with Default}
\author{N. Agram$^{1,2}$, S. Labed$^{2}$, B. Mansouri$^{2}$ \& M. A. Saouli$^{2}$}
\maketitle

\begin{abstract}
We are interested on reflected advanced backward stochastic differential
equations (RABSDE) with default. By the predictable representation property
and for a Lipschitz driver, we show that the RABSDE with default has a unique
solution in the enlarged filtration. A comparison theorem for such type of
equations is proved. Finally, we give a connection between RABSDE and optimal stopping.

\end{abstract}


\footnotetext[1]{Department of Mathematics, University of Oslo, P.O. Box 1053
Blindern, N--0316 Oslo, Norway.\newline Email: \texttt{naciraa@math.uio.no.}
\par
This research was carried out with support of the Norwegian Research Council,
within the research project Challenges in Stochastic Control, Information and
Applications (STOCONINF), project number 250768/F20.
}

\footnotetext[2]{Department of Mathematics, University of Biskra, Algeria.
\par
Emails: labed.saloua@yahoo.fr, mansouri.badreddine@gmail.com,
saoulimoustapha@yahoo.fr.}


\section{Introduction}

Reflected advanced backward stochastic differential equations (RABSDE) appear
in their linear form as the adjoint equation when dealing with the stochastic
maximum principle to study optimal singular control for delayed systems, we
refer for example to Øksendal and Sulem \cite{os-sing} and also to Agram
\textit{et al} \cite{abop} for more general case. This is a natural model in
population growth, but also in finance, where people's memory plays a role in
the price dynamics.\newline

After the economic crises in $2008$, researchers started to include default in
banks as a part of their financial modelling. This is why we are interested on
RABSDE also in the context of enlargement of filtration. In order to be more
precise, let us consider a random time $\tau$ which is neither an $\mathbb{F}%
$-stopping time nor $\mathcal{F}_{T}$-measurable. Examples of such random
times are default times, where the reason for the default comes from outside
the Brownian model. We denote $H_{t}=\mathbf{1}_{\tau\leq t},\ t\in
\lbrack0,T],$ and consider the filtration $\mathbb{G}$ obtained by enlarging
progressively the filtration $\mathbb{F}$ by the process $H$, i.e.,
$\mathbb{G}$ is the smallest filtration satisfying the usual assumptions of
completeness and right-continuity, which contains the filtration $\mathbb{F}$
and has $H$ as an adapted process. The RABSDE related with, we want to study
is the following:%

\[
\left\{
\begin{array}
[c]{ll}%
Y_{t} & =\xi+%
{\textstyle\int_{t}^{T}}
f(s,Y_{s},Z_{s},\mathbb{E[}Y_{s+\delta}|\mathcal{G}_{s}],\mathbb{E[}%
Z_{s+\delta}|\mathcal{G}_{s}],U_{s})ds-%
{\textstyle\int_{t}^{T}}
Z_{s}dW_{s}\\
& -%
{\textstyle\int_{t}^{T}}
U_{s}dH_{s}+K_{T}-K_{t},\qquad t\in\left[  0,T\right]  ,\\
Y_{t} & =\xi,\qquad t\geq T,\\
Z_{t} & =U_{t}=0,\qquad t>T.
\end{array}
\right.
\]

By saying that the RBSDE is advanced we mean that driver at the present time
$s$ may depend not only on present values of the solution processes $(Y,Z,K)$,
but also on the future values $s+\delta$ for some $\delta> 0$. To make the
system adapted, we take the conditional expectation of the advanced terms.

We will see that by using the predictable representation property (PRP) the
above system is equivalent to a RABSDE driven by a martingale, consisting of
the Brownian motion $W$ and the martingale $M$ associated to the jump process
$H$, as follows:
\[
\left\{
\begin{array}
[c]{ll}%
Y_{t} & =\xi+%
{\textstyle\int_{t}^{T}}
F(s,Y_{s},Z_{s},\mathbb{E[}Y_{s+\delta}|\mathcal{G}_{s}],\mathbb{E[}%
Z_{s+\delta}|\mathcal{G}_{s}],U_{s})ds-%
{\textstyle\int_{t}^{T}}
Z_{s}dW_{s}\\
& -%
{\textstyle\int_{t}^{T}}
U_{s}dM_{s}+K_{T}-K_{t},\qquad t\in\left[  0,T\right]  ,\\
Y_{t} & =\xi,\qquad t\geq T,\\
Z_{t} & =U_{t}=0,\qquad t>T.
\end{array}
\right.
\]

Our aim in this paper is not to find solutions in the Brownian filtration by
using the decomposition approach, as it has been done in Kharroubi and Lim
\cite{elkarui} for BSDE and in Jeanblanc \textit{et al }\cite{jla} for ABSDE.
However, we want to find solutions under the enlarged filtration rather than
the Brownian one as in the previous works.

In Dumitrescu \textit{et al }\cite{dqs1}, \cite{dqs2}, \cite{dqs3}, the
authors consider directly BSDE and RBSDE driven by general filtration
generated by the pair $(W,M)$.

We will extend the recent woks by Dumitrescu \textit{et al }\cite{dqs1},
\cite{dqs2}, \cite{dqs3}\textit{, }to the anticipated case and we will explain
how such an equations appear by using the PRP.

We will extend also the comparison theorem for ABSDE in Peng and Yang
\cite{peng2009} to RABSDE with default and finally, we give a link between
RABSDE with default and optimal stopping as it has been done in El Karoui
\textit{et al }\cite{elkarui} and Øksendal and Zhang \cite{oz}.

For more details about ABSDE with jumps coming from the compensated Poisson
random measure which is independent of the Brownian motion, we refer to
Øksendal \textit{et al} \cite{Oksendal2011}, \cite{os}. For RBSDE with jumps,
we refer to Quenez and Sulem \cite{qs} and for more details about enlargement
progressive of filtration, we refer to Song \cite{songspliting} .

\section{Framework}

Let $(\Omega,\mathcal{G},P)$ be a complete probability space. We assume that
this space is equipped with a one-dimensional standard Brownian motion $W$ and
we denote by $\mathbb{F}:=(\mathcal{F}_{t})_{t\geq0}$ the right continuous
complete filtration generated by $W$. We also consider on this space a random
time $\tau$, which represents for example a default time in credit risk or in
counterparty risk, or a death time in actuarial issues. The random time $\tau$
is not assumed to be an $\mathbb{F}$-stopping time. We therefore use in the
sequel the standard approach of filtration enlargement by considering
$\mathbb{G}$ the smallest right continuous extension of $\mathbb{F}$ that
turns $\tau$ into a $\mathbb{G}$-stopping time (see e.g. Chapter 4 in
\cite{aj}). More precisely $\mathbb{G}:=(\mathcal{G}_{t})_{t\geq0}$ is defined
by
\[
\mathcal{G}_{t}:=\bigcap_{\varepsilon>0}\tilde{\mathcal{G}}_{t+\varepsilon
}\;,
\]
for all $t\geq0$, where $\tilde{\mathcal{G}}_{s}:=\mathcal{F}_{s}\vee
\sigma(\mathbf{1}_{\tau\leq u}\;,u\in\lbrack0,s])$, for all $s\geq0$.

\vspace{2mm}

We denote by $\mathcal{P}(\mathbb{G})$ the $\sigma$-algebra of $\mathbb{G}%
$-predictable subsets of $\Omega\times\lbrack0,T]$, i.e. the $\sigma$-algebra
generated by the left-continuous $\mathbb{G}$-adapted processes.

\vspace{2mm}

We then impose the following assumptions, which are classical in the
filtration enlargement theory.

\begin{description}
\item[$(H)$] The process $W$ is a $\mathbb{G}$-Brownian motion. We observe
that, since the filtration $\mathbb{F}$ is generated by the Brownian motion
$W$, this is equivalent with the fact that all $\mathbb{F}$-martingales are
also $\mathbb{G}$-martingales. Moreover, it also follows that the stochastic
integral $%
{\textstyle\int_{0}^{t}}
X_{s}dW_{s}$ is well defined for all $\mathcal{P}(\mathbb{G})$-measurable
processes $X$ such that $%
{\textstyle\int_{0}^{t}}
\left\vert X_{s}\right\vert ^{2}ds<\infty,$ for all $t\geq0$.
\end{description}

\begin{itemize}
\item The process $M$ defined by
\[
M_{t}=H_{t}-%
{\textstyle\int_{0}^{t\wedge\tau}}
\lambda_{s}ds,\qquad t\geq0,
\]
is a $\mathbb{G}$-martingale with single jump time $\tau$ and the process
$\lambda$ is ${\mathbb{F}}$-adapted, called the ${\mathbb{F}}$-intensity of
$\tau$. \newline

\item We assume that the process $\lambda$ is upper bounded by a constant.

\item Under $(H)$ any square integrable ${\mathbb{G}}$ martingale $Y$ admits a
representation as
\[
Y_{t}=y+%
{\textstyle\int_{0}^{t}}
\varphi_{s}dW_{s}+%
{\textstyle\int_{0}^{t}}
\gamma_{s}dM_{s},
\]
where $M$ is the compensated martingale of $H$, and $\varphi,\gamma$ are
square-integrable ${\mathbb{G}}$-predictable processes. (See Theorem 3.15 in
\cite{aj}).\vspace{2mm}
\end{itemize}

Throughout this section, we introduce some basic notations and spaces.

\begin{itemize}
\item $S_{\mathbb{G}}^{2}$ is the subset of $%
\mathbb{R}
$-valued $\mathbb{G}$-adapted càdlàg processes $\left(  Y_{t}\right)
_{t\in\left[  0,T\right]  }$, such that
\[
\left\Vert Y\right\Vert _{S^{2}}^{2}:=\mathbb{E[}\underset{t\in\left[
0,T\right]  }{\sup}\left\vert Y_{t}\right\vert ^{2}]<\infty.
\]

\item $\mathcal{K}^{2}$ is a set of real-valued nondecreasing processes $K$
with $K_{0^{-}}=0$ and $\mathbb{E[}K_{t}]<\infty.$

\item $H_{\mathbb{G}}^{2}$ is the subset of $%
\mathbb{R}
$-valued $\mathcal{P}(\mathbb{G})$-measurable processes $\left(  Z_{t}\right)
_{t\in\left[  0,T\right]  },$ such that
\[
\left\Vert Z\right\Vert _{H^{2}}^{2}:=\mathbb{E[}%
{\textstyle\int_{0}^{T}}
\left\vert Z_{t}\right\vert ^{2}dt]<\infty.
\]

\item $L^{2}(\lambda)$ is the subset of $%
\mathbb{R}
$-valued $\mathcal{P}(\mathbb{G})$-measurable processes $\left(  U_{t}\right)
_{t\in\left[  0,T\right]  },$ such that
\[
\left\Vert U\right\Vert _{L^{2}(\lambda)}^{2}:=\mathbb{E[}%
{\textstyle\int_{0}^{T\wedge\tau}}
\lambda_{t}\left\vert U_{t}\right\vert ^{2}dt]<\infty.
\]

\end{itemize}

\section{Existence and Uniqueness}

We study the RABSDE with default%

\begin{equation}
\left\{
\begin{array}
[c]{ll}%
Y_{t} & =\xi+%
{\textstyle\int_{t}^{T}}
f(s,Y_{s},Z_{s},\mathbb{E[}Y_{s+\delta}|\mathcal{G}_{s}],\mathbb{E[}%
Z_{s+\delta}|\mathcal{G}_{s}],U_{s})ds-%
{\textstyle\int_{t}^{T}}
Z_{s}dW_{s}\\
& -%
{\textstyle\int_{t}^{T}}
U_{s}dH_{s}+K_{T}-K_{t},\qquad t\in\left[  0,T\right]  ,\\
Y_{t} & =\xi,\qquad t\geq T,\\
Z_{t} & =U_{t}=0,\qquad t>T,
\end{array}
\right.  \label{4}%
\end{equation}
where $f$ is $\mathcal{G}_{t}\otimes\mathcal{B}\left(  \left[  0,T\right]
\right)  \otimes\mathcal{B}\left(
\mathbb{R}
^{5}\right)  $-measurable, and the terminal condition $\xi$ is $\mathcal{G}%
_{T}$-measurable. Moreover

\begin{itemize}
\item $Y_{t}\geq S_{t}$, for each $t\geq0$ a.s.

\item $K_{t}$ is càdlàg, increasing and $\mathbb{G}$-adapted process with
$K_{0^{-}}=0.$

\item $%
{\textstyle\int_{0}^{T}}
(Y_{t}-S_{t})dK_{t}^{c}=0$ and $\triangle K_{t}^{d}=-\triangle Y_{t}%
\mathbf{1}_{\{Y_{t^{-}}=S_{t^{-}}\}}$, where denote the continuous and
discontinuous parts of $K$ respectively.

\item $(S_{t})_{t\geq0}$ is the obstacle which is a càdlàg, increasing and
$\mathbb{G}$-adapted process.
\end{itemize}

We call the quadruplet $(Y,Z,U,K)$ solution of the RABSDE (\ref{4}%
).\vspace{2mm}

Let us impose the following set of assumptions.

(\textit{i) Assumption on the terminal condition:}

\begin{itemize}
\item $\xi\in L^{2}\left(  \Omega,\mathcal{G}_{T}\right)  $.
\end{itemize}

(ii) \textit{Assumptions on the generator function }$f:\Omega\times\left[
0,T\right]  \times\mathcal{\mathbb{R}}^{5}\rightarrow\mathbb{R}$ is such that

\begin{itemize}
\item $\mathbb{G}$-predictable and satisfies the integrability condition, such
that
\begin{equation}
\mathbb{E[}%
{\textstyle\int_{0}^{T}}
|f(t,0,0,0,0,0)|^{2}dt]<0\text{, } \label{int}%
\end{equation}

\end{itemize}

for all $t\in\left[  0,T\right]  .$

\begin{itemize}
\item Lipschitz in the sense that, there exists $C>0,$ such that%
\begin{align}
&  |f(t,y,z,\mu,\pi,u)-f(t,y^{\prime},z^{\prime},\mu^{\prime},\pi^{\prime
},u^{\prime})|\nonumber\\
&  \leq C(\left\vert y-y^{\prime}\right\vert +|z-z^{\prime}|+\left\vert
\pi-\pi^{\prime}\right\vert +\left\vert \mu-\mu^{\prime}\right\vert
+\lambda_{t}|u-u^{\prime}|), \label{lip}%
\end{align}

\end{itemize}

for all $t\in\left[  0,T\right]  $ and all $y,y^{\prime},z,z^{\prime},\mu
,\mu^{\prime},\pi,\pi^{\prime},u,u^{\prime}\in%
\mathbb{R}
.$\vspace{2mm}

We give the existence of the solution to a RABSDE in the enlarged filtration
$\mathbb{G}.$ The existence follows from the PRP as we can also say, the
property of martingale representation (PMR), and a standard approach like any
classical RBSDE.

\vspace{2mm}

Under our assumptions we know that equation (\ref{4}) is equivalent to%
\begin{equation}
\left\{
\begin{array}
[c]{ll}%
Y_{t} & =\xi+%
{\textstyle\int_{t}^{T}}
F(s,Y_{s},Z_{s},\mathbb{E[}Y_{s+\delta}|\mathcal{G}_{s}],\mathbb{E[}%
Z_{s+\delta}|\mathcal{G}_{s}],U_{s})ds-%
{\textstyle\int_{t}^{T}}
Z_{s}dW_{s}\\
& -%
{\textstyle\int_{t}^{T}}
U_{s}dM_{s}+K_{T}-K_{t},\qquad t\in\left[  0,T\right]  ,\\
Y_{t} & =\xi,\qquad t\geq T,\\
Z_{t} & =U_{t}=0,\qquad t>T,
\end{array}
\right.  \label{5}%
\end{equation}
with $dH_{s}=dM_{s}+\lambda_{s}\mathbf{1}_{s<\tau}ds$, and
\[
F(s,y,z,\mu,\pi,u):=f(s,y,z,\mu,\pi^{\prime},u)-\lambda_{s}(\mathbf{1}%
-H_{s})u.
\]
By assumption, the process $\lambda$ is bounded.

\vspace{2mm}

In order to get existence and uniqueness for the RABSDE (\ref{5}), let us
check that the generator $F$ satisfies the same assumption as $f:$ The
function\textit{ }$F:\Omega\times\left[  0,T\right]  \times\mathcal{\mathbb{R}%
}^{5}\rightarrow\mathbb{R}$ is such that

\begin{description}
\item[(i)] $\mathbb{G}$-predictable and integrable in the sense that, for all
$t\in\left[  0,T\right]  $, by inequality (\ref{int}), we have
\[
\mathbb{E[}%
{\textstyle\int_{0}^{T}}
|F(t,0,0,0,0,0)|^{2}dt]=\mathbb{E[}%
{\textstyle\int_{0}^{T}}
|f(t,0,0,0,0,0)|^{2}dt]<0.
\]

\item[(ii)] Lipschitz in the sense that there exists a constant $C^{\prime}%
>0$, such that%
\begin{equation}%
\begin{array}
[c]{l}%
|F(t,y,z,\mu,\pi,u)-F(t,y^{\prime},z^{\prime},\mu^{\prime},\pi^{\prime
},u^{\prime})|\\
=|f(t,y,z,\mu,\pi,u)-f(t,y^{\prime},z^{\prime},\mu^{\prime},\pi^{\prime
},u^{\prime})-\lambda_{t}(\mathbf{1}-H_{t})(u-u^{\prime})|\\
\leq|f(t,y,z,\mu,\pi,u)-f(t,y^{\prime},z^{\prime},\mu^{\prime},\pi^{\prime
},u^{\prime})|+\lambda_{t}(\mathbf{1}-H_{t})|u-u^{\prime}|\\
\leq C(\left\vert y-y^{\prime}\right\vert +|z-z^{\prime}|+\left\vert \mu
-\mu^{\prime}\right\vert +\left\vert \pi-\pi^{\prime}\right\vert +\lambda
_{t}(\mathbf{1}-H_{t})|u-u^{\prime}|)+\lambda_{t}(\mathbf{1}-H_{t}%
)|u-u^{\prime}|\\
\leq C^{\prime}(\left\vert y-y^{\prime}\right\vert +|z-z^{\prime}|+\left\vert
\mu-\mu^{\prime}\right\vert +\left\vert \pi-\pi^{\prime}\right\vert
+\lambda_{t}|u-u^{\prime}|),
\end{array}
\nonumber
\end{equation}

\end{description}

for all $t\in\left[  0,T\right]  $ and all $y,z,u,,\pi,\mu,y^{\prime
},z^{\prime},u^{\prime},\pi^{\prime},\mu^{\prime}\in%
\mathbb{R}
,$where we have used the Lipschitzianity of $f$ (\ref{lip}).

\begin{description}
\item[(iii)] The terminal value: $\xi\in L^{2}\left(  \Omega,\mathcal{G}%
_{T}\right)  $.\bigskip
\end{description}

\begin{theorem}
\label{ex-uni}Under the above assumptions (i)-(iii), the RABSDE $\left(
\ref{5}\right)  $ admits a unique solution $\left(  Y,Z,U,K\right)  \in
S_{\mathbb{G}}^{2}\times H_{\mathbb{G}}^{2}\times L^{2}(\lambda)\times
\mathcal{K}^{2}.$
\end{theorem}

\noindent{Proof.} \quad We define the mapping\textbf{\ }%
\[
\Phi:H_{\mathbb{G}}^{2}\times H_{\mathbb{G}}^{2}\times L^{2}(\lambda
)\rightarrow H_{\mathbb{G}}^{2}\times H_{\mathbb{G}}^{2}\times L^{2}%
(\lambda),
\]
for which we will show that it is contracting under a suitable norm. For this
we note that for any $\left(  Y,Z,U,K\right)  \in H_{\mathbb{G}}^{2}\times
H_{\mathbb{G}}^{2}\times L^{2}(\lambda)\times\mathcal{K}^{2}$ there exists a
unique quadruple $(\hat{Y},\hat{Z},\hat{U},\hat{K})\in S_{\mathbb{G}}%
^{2}\times H_{\mathbb{G}}^{2}\times L^{2}(\lambda)\times\mathcal{K}^{2},$ such
that%
\begin{equation}%
\begin{array}
[c]{c}%
\hat{Y}_{t}=\xi+%
{\textstyle\int_{t}^{T}}
F(s,Y_{s},Z_{s},\mathbb{E[}Y_{s+\delta}|\mathcal{G}_{s}],\mathbb{E[}%
Z_{s+\delta}|\mathcal{G}_{s}],U_{s})ds-%
{\textstyle\int_{t}^{T}}
\hat{Z}_{s}dW_{s}\\
-%
{\textstyle\int_{t}^{T}}
\hat{U}_{s}dM_{s}-\int_{t}^{T}d\hat{K}_{s},\qquad t\in\left[  0,T\right]  ,
\end{array}
\label{n1}%
\end{equation}
Let $\Phi(Y,Z,U):=(\hat{Y},\hat{Z},\hat{U})$. For given $\left(  Y^{i}%
,Z^{i},U^{i}\right)  $\textbf{\ }$\in H_{\mathbb{F}}^{2}\times H_{\mathbb{F}%
}^{2}\times L^{2}(\lambda),$\textbf{\ }for $i=1,2$, we use the simplified
notations:
\[%
\begin{array}
[c]{ll}%
(\hat{Y}^{i},\hat{Z}^{i},\hat{U}^{i}) & :=\Phi(Y^{i},Z^{i},U^{i}),\\
(\tilde{Y},\tilde{Z},\tilde{U}) & :=(\hat{Y}^{1},\hat{Z}^{1},\hat{U}%
^{1})-(\hat{Y}^{2},\hat{Z}^{2},\hat{U}^{2}),\\
(\bar{Y},\bar{Z},\bar{U}) & :=(Y^{1},Z^{1},U^{1})-(Y^{2},Z^{2},U^{2}).
\end{array}
\]
The triplet of processes $\left(  \tilde{Y},\tilde{Z},\tilde{U}\right)  $
satisfies the equation%
\[%
\begin{array}
[c]{c}%
\tilde{Y}_{t}=%
{\textstyle\int_{t}^{T}}
\mathbb{\{}F(s,Y_{s}^{1},Z_{s}^{1},\mathbb{E[}Y_{s+\delta}^{1}|\mathcal{G}%
_{s}],\mathbb{E[}Z_{s+\delta}^{1}|\mathcal{G}_{s}],U_{s}^{1})\\
-F(s,Y_{s}^{2},Z_{s}^{2},\mathbb{E[}Y_{s+\delta}^{2}|\mathcal{G}%
_{s}],\mathbb{E[}Z_{s+\delta}^{2}|\mathcal{G}_{s}],U_{s}^{2})\}ds\\
-%
{\textstyle\int_{t}^{T}}
\tilde{Z}_{s}dW_{s}-%
{\textstyle\int_{t}^{T}}
\tilde{U}_{s}dM_{s}-\int_{t}^{T}d\tilde{K}_{s},\qquad t\in\left[  0,T\right]
.
\end{array}
\]
We have that $M_{t}=H_{t}-%
{\textstyle\int_{0}^{t}}
\lambda_{s}ds$ which is a pure jump martingale. Then,%
\[
\left[  M\right]  _{t}=\underset{0\leq s\leq t}{%
{\textstyle\sum}
}\left(  \bigtriangleup M_{s}\right)  ^{2}=\underset{0\leq s\leq t}{%
{\textstyle\sum}
}\left(  \bigtriangleup H_{s}\right)  ^{2}=H_{t},
\]
and%
\[
\left\langle M\right\rangle _{t}=%
{\textstyle\int_{0}^{t}}
\lambda_{s}ds,
\]%
\[%
{\textstyle\int_{t}^{T}}
|\tilde{U}_{s}|^{2}d\left\langle M\right\rangle _{s}=%
{\textstyle\int_{t}^{T}}
\lambda_{s}|\tilde{U}_{s}|^{2}ds.
\]
Applying Itô's formula to\textbf{\ }$e^{\beta t}|\tilde{Y}_{t}|^{2}$\textbf{,}
taking conditional expectation and using the Lipschitz condition, we get%

\[%
\begin{array}
[c]{l}%
\mathbb{E[}%
{\textstyle\int_{0}^{T}}
e^{\beta s}(\beta|\tilde{Y}_{s}|^{2}+|\tilde{Z}_{s}|^{2}+\lambda_{s}|\tilde
{U}_{s}|^{2})ds]\\
\leq10\rho C^{2}\mathbb{E[}%
{\textstyle\int_{0}^{T}}
e^{\beta s}\left\vert \bar{Y}_{s}\right\vert ^{2}ds]+\tfrac{1}{2\rho
}\mathbb{E[}%
{\textstyle\int_{0}^{T}}
e^{\beta s}\{\left\vert \bar{Z}_{s}\right\vert ^{2}+\lambda_{s}^{2}\left\vert
\bar{U}_{s}\right\vert ^{2}\}ds],
\end{array}
\]
where we have used that
\begin{align*}
\tilde{Y}_{s}dK_{s}^{1,c}  &  =(Y_{s}^{1}-S_{s})dK_{s}^{1,c}-(Y_{s}^{2}%
-S_{s})dK_{s}^{1,c}\\
&  =-(Y_{s}^{2}-S_{s})dK_{s}^{1,c}\leq0\text{ a.s.,}%
\end{align*}
and by symmetry, we have also $\tilde{Y}_{s}dK_{s}^{2,c}\geq0$ a.s. For the
discontinuous case, we have as well%
\begin{align*}
\tilde{Y}_{s}dK_{s}^{1,d}  &  =(Y_{s}^{1}-S_{s})dK_{s}^{1,d}-(Y_{s}^{2}%
-S_{s})dK_{s}^{1,d}\\
&  =-(Y_{s}^{2}-S_{s})dK_{s}^{1,d}\leq0\text{ a.s.,}%
\end{align*}
and by symmetry, we have also $\tilde{Y}_{s}dK_{s}^{2,d}\geq0$ a.s.

\noindent Since $\lambda$ is bounded, we get that $\lambda^{2}\leq k\lambda$
and by choosing $\beta=1+10\rho C^{2}\ $we obtain%
\[
||(\tilde{Y},\tilde{Z},\tilde{U})||^{2}\leq\tfrac{1}{2\rho}||(\bar{Y},\bar
{Z},\bar{U})||^{2}%
\]
which means for $\rho\geq1$, there exists a unique fixed point that is a
solution for our RABSDE $\left(  \ref{5}\right)  .\qquad\square$ \newline

\section{Comparison Theorem for RABSDE with Default}

In this section we are interested in a subclass of RABSDE where the driver
only depend on future values of $Y$ and is not allowed to depend on future
values of $Z$, as follows:%
\[
\left\{
\begin{array}
[c]{ll}%
Y_{t} & =\xi+%
{\textstyle\int_{t}^{T}}
g(s,Y_{s},Z_{s},\mathbb{E[}Y_{s+\delta}|\mathcal{G}_{s}],U_{s})ds-%
{\textstyle\int_{t}^{T}}
Z_{s}dW_{s}\\
& -%
{\textstyle\int_{t}^{T}}
U_{s}dM_{s}+K_{T}-K_{t},\qquad t\in\left[  0,T\right]  ,\\
Y_{t} & =\xi,\qquad t\geq T,\\
Z_{t} & =U_{t}=0,\qquad t>T,
\end{array}
\right.
\]
such that

\begin{itemize}
\item $Y_{t}\geq S_{t}$, for each $t\geq0$ a.s.

\item $K_{t}$ is càdlàg, increasing and $\mathbb{G}$-adapted process with
$K_{0^{-}}=0.$

\item $%
{\textstyle\int_{0}^{T}}
(Y_{t}-S_{t})dK_{t}^{c}=0$ and $\triangle K_{t}^{d}=-\triangle Y_{t}%
\mathbf{1}_{\{Y_{t^{-}}=S_{t^{-}}\}}$, where denote the continuous and
discontinuous parts of $K$ respectively.

\item $(S_{t})_{t\geq0}$ is the obstacle which is a càdlàg, increasing and
$\mathbb{G}$-adapted process.\newline\newline
\end{itemize}

We impose the following set of assumptions.

\begin{description}
\item[(a)] The driver $g:\Omega\times\left[  0,T\right]  \times
\mathcal{\mathbb{R}}^{4}\rightarrow\mathbb{R}$ is $\mathbb{G}$-predictable,
and satisfies
\[
\mathbb{E[}%
{\textstyle\int_{0}^{T}}
|g(t,0,0,0,0)|^{2}dt]<0\text{, }%
\]

\end{description}

\begin{align*}
&  |g(t,y,z,\mu,u)-g(t,y^{\prime},z^{\prime},\mu^{\prime},u^{\prime})|\\
&  \leq C(\left\vert y-y^{\prime}\right\vert +|z-z^{\prime}|+\left\vert
\mu-\mu^{\prime}\right\vert +\lambda_{t}|u-u^{\prime}|),
\end{align*}

for all $t\in\left[  0,T\right]  $ and all $y,y^{\prime},z,z^{\prime},\mu
,\mu^{\prime},u,u^{\prime}\in%
\mathbb{R}
.$

\begin{description}
\item[(b)] T\textit{he terminal condition: }$\xi\in L^{2}\left(
\Omega,\mathcal{G}_{T}\right)  $.\newline\newline
\end{description}

Let us first state the comparison theorem for RBSDE with default which relies
on the comparison theorem for BSDE with default done by Dumitrescu \textit{et
al} \cite{dqs2}, Theorem 2.17.

\begin{theorem}
\label{comp} Let $g^{1},g^{2}:\Omega\times\left[  0,T\right]  \times
\mathcal{\mathbb{R}}^{3}\rightarrow\mathbb{R},$ $\xi^{1},\xi^{2}\in
L^{2}\left(  \Omega,\mathcal{G}_{T}\right)  $ and let the quadruplet
$(Y^{j},Z^{j},U^{j},K^{j})_{j=1,2}$ be the solution of the RBSDE with default%
\[
\left\{
\begin{array}
[c]{ll}%
Y_{t}^{j} & =\xi^{j}+\int_{t}^{T}g^{j}(s,Y_{s}^{j},Z_{s}^{j},U_{s}^{j})ds-%
{\textstyle\int_{t}^{T}}
Z_{s}^{j}dW_{s}\\
& -%
{\textstyle\int_{t}^{T}}
U_{s}^{j}dM_{s}+\int_{t}^{T}dK_{s}^{j},\ \qquad t\in\left[  0,T\right]  ,\\
Y_{t}^{j} & =\xi^{j},\ \qquad t\geq T,\\
Z_{t}^{j} & =U_{t}^{j}=0,\ \qquad t>T.
\end{array}
\right.
\]
The drivers $(g^{j})_{j=1,2}$ satisfies assumptions (a)-(b). Suppose that
there exists a predictable process $(\theta_{t})_{t\geq0}$ with $\theta
_{t}\lambda_{t}$ bounded and $\theta_{t}\geq-1$ $dt\otimes dP$ a.s. such that%
\[
g^{1}\left(  t,y,z,u\right)  -g^{1}\left(  t,y,z,u^{\prime}\right)  \geq
\theta_{t}(u-u^{\prime})\lambda_{t}.
\]
Moreover, suppose that
\end{theorem}

\begin{itemize}
\item $\xi^{1}\geq\xi^{2}$, a.s.

\item For any\textbf{\ }$t\in\left[  0,T\right]  ,$ $S_{t}^{1}\geq S_{t}^{2}$, a.s.

\item $g^{1}\left(  t,y,z,u\right)  \geq g^{2}\left(  t,y,z,u\right)  ,$
$\forall t\in\left[  0,T\right]  ,$ $y,z,u\in%
\mathbb{R}
.$\newline Then
\[
Y_{t}^{1}\geq Y_{t}^{2},\qquad\forall t\in\left[  0,T\right]  \text{.}%
\]

\end{itemize}

\begin{theorem}
\label{theorem 4.1} Let $g^{1},g^{2}:\Omega\times\left[  0,T\right]
\times\mathcal{\mathbb{R}}^{4}\rightarrow\mathbb{R},$ $\xi^{1},\xi^{2}\in
L^{2}\left(  \Omega,\mathcal{G}_{T}\right)  $ and let the quadruplet
$(Y^{j},Z^{j},U^{j},K^{j})_{j=1,2}$ be the solution of the RABSDE%
\[
\left\{
\begin{array}
[c]{ll}%
Y_{t}^{j} & =\xi^{j}+\int_{t}^{T}g^{j}(s,Y_{s}^{j},Z_{s}^{j},\mathbb{E[}%
Y_{s+\delta}^{j}|\mathcal{G}_{s}],U_{s}^{j})ds-%
{\textstyle\int_{t}^{T}}
Z_{s}^{j}dW_{s}\\
& -%
{\textstyle\int_{t}^{T}}
U_{s}^{j}dM_{s}+\int_{t}^{T}dK_{s}^{j},\ \qquad t\in\left[  0,T\right]  ,\\
Y_{t}^{j} & =\xi^{j},\ \qquad t\geq T,\\
Z_{t}^{j} & =U_{t}^{j}=0,\ \qquad t>T.
\end{array}
\right.
\]
The drivers $(g^{j})_{j=1,2}$ satisfies assumptions (a)-(b). Moreover, suppose that:
\end{theorem}

\begin{itemize}
\item[(i)] For all $t\in\left[  0,T\right]  ,$ $y,$ $z,u\in%
\mathbb{R}
,$ $g^{2}\left(  t,y,z,\cdot,u\right)  $ is increasing with respect to
$Y_{t+\delta}$ in the sense that%
\[
g^{2}\left(  t,y,z,Y_{t+\delta},u\right)  \geq g^{2}\left(  t,y,z,Y_{t+\delta
}^{\prime},u\right)  ,
\]
for all $Y_{t+\delta}\geq Y_{t+\delta}^{\prime}.$

\item[(ii)] $\xi^{1}\geq\xi^{2}$, $a.s$.

\item[(iii)] For each\textbf{\ }$t\in\left[  0,T\right]  ,$ $S_{t}^{1}\geq
S_{t}^{2}$, a.s.

\item[(iv)] Suppose that there exists a predictable process $(\theta
_{t})_{t\geq0}$ with $\theta_{t}\lambda_{t}$ bounded and $\theta_{t}\geq-1$
$dt\otimes dP$ a.s., such that%
\[
g^{1}\left(  t,y,z,Y_{t+\delta},u\right)  -g^{1}\left(  t,y,z,Y_{t+\delta
},u^{\prime}\right)  \geq\theta_{t}(u-u^{\prime})\lambda_{t}.
\]

\item[(v)] $g^{1}\left(  t,y,z,Y_{t+\delta},u\right)  \geq g^{2}\left(
t,y,z,Y_{t+\delta},u\right)  ,$ $\forall t\in\left[  0,T\right]  ,$
$y,z,Y_{t+\delta},u\in%
\mathbb{R}
.$\newline Then, we have
\[
Y_{t}^{1}\geq Y_{t}^{2},\qquad\text{a.e.,a.s.}%
\]

\end{itemize}

\noindent{Proof.} \quad\ Consider the following RABSDE%

\[
\left\{
\begin{array}
[c]{ll}%
Y_{t}^{3} & =\xi^{2}+\int_{t}^{T}g^{2}(s,Y_{s}^{3},Z_{s}^{3},\mathbb{E[}%
Y_{s+\delta}^{1}|\mathcal{G}_{s}],U_{s}^{3})ds-%
{\textstyle\int_{t}^{T}}
Z_{s}^{3}dW_{s}\\
& -%
{\textstyle\int_{t}^{T}}
U_{s}^{3}dM_{s}+\int_{t}^{T}dK_{s}^{3},\text{\ \qquad\qquad}t\in\left[
0,T\right]  ,\\
Y_{t}^{3} & =\xi^{2},\text{\ \qquad\qquad}t\geq T,\\
Z_{t}^{3} & =U_{t}^{3}=0,\text{\ \qquad\qquad}t>T.
\end{array}
\right.
\]
From Proposition 3.2 in Dumitrescu \textit{et al} \cite{dqs3}, we know there
exists a unique quadruplet of $\mathbb{G}$-adapted processes $\left(
Y^{3},Z^{3},U^{3},K^{3}\right)  \in S_{\mathbb{G}}^{2}\times H_{\mathbb{G}%
}^{2}\times L^{2}(\lambda)\times\mathcal{K}^{2}$ satisfies the above RBSDE
since the advanced term is considered as a parameter.

\noindent Now we have by assumptions (iii)-(v) and Theorem \ref{comp}, that
\[
Y_{t}^{1}\geq Y_{t}^{3},\qquad\text{ for all }t\text{, a.s.}%
\]
Set
\[
\left\{
\begin{array}
[c]{ll}%
Y_{t}^{4} & =\xi^{2}+%
{\textstyle\int_{t}^{T}}
g^{2}(s,Y_{s}^{4},Z_{s}^{4},\mathbb{E}[Y_{s+\delta}^{3}|\mathcal{G}_{s}%
],U_{s}^{4})ds-%
{\textstyle\int_{t}^{T}}
Z_{s}^{4}dW_{s}\\
& \qquad-%
{\textstyle\int_{t}^{T}}
U_{s}^{4}dM_{s}+\int_{t}^{T}dK_{s}^{4},\text{\ \qquad}t\in\left[  0,T\right]
,\\
Y_{t}^{4} & =\xi^{2},\text{\qquad}t\geq T,\\
Z_{t}^{4} & =U_{t}^{4}=0,\text{\qquad}t>T.
\end{array}
\right.
\]
By the same arguments, we get $Y_{t}^{3}\geq Y_{t}^{4},$ a.e.,a.s.

\noindent For $n=5,6,...,$ we consider the following RABSDE%
\[
\left\{
\begin{array}
[c]{ll}%
Y_{t}^{n} & =\xi^{2}+%
{\textstyle\int_{t}^{T}}
g^{2}(s,Y_{s}^{n},Z_{s}^{n},\mathbb{E[}Y_{s+\delta}^{n-1}|\mathcal{G}%
_{t}],U_{s}^{n})ds-%
{\textstyle\int_{t}^{T}}
Z_{s}^{n}dW_{s}\\
& -%
{\textstyle\int_{t}^{T}}
U_{s}^{n}dM_{s}+%
{\textstyle\int_{t}^{T}}
dK_{s}^{n},\text{\ \qquad}t\in\left[  0,T\right]  ,\\
Y_{t}^{n} & =\xi^{2},\text{\qquad}t\geq T,\\
Z_{t}^{n} & =U_{t}^{n}=0,\ \qquad t>T.
\end{array}
\right.
\]
We may remark that it is clear that $Y_{s+\delta}^{n-1}$ is considered to be
knowing on the above RABSDE. By induction on $n>4$, we get%
\[
Y_{t}^{4}\geq Y_{t}^{5}\geq Y_{t}^{6}\geq\cdot\cdot\cdot\geq Y_{t}^{n}%
\geq\cdot\cdot\cdot,\text{\quad a.s.}%
\]
If we denote by%
\[
\bar{Y}=Y^{n}-Y^{n-1},\text{\qquad}\bar{Z}=Z^{n}-Z^{n-1},\text{\qquad}\bar
{U}=U^{n}-U^{n-1},\text{\qquad}\bar{K}=K^{n}-K^{n-1}.
\]
By similar estimations as in the proof of Theorem \ref{ex-uni}, we can find
that $\left(  Y^{n},Z^{n},U^{n},K^{n}\right)  $ converges to $\left(
Y^{n-1},Z^{n-1},U^{n-1},K^{n-1}\right)  $ as $n\rightarrow\infty.$

\noindent Iterating with respect to $n$, we obtain when $n\rightarrow\infty,$
that $\left(  Y^{n},Z^{n},U^{n},K^{n}\right)  $ converges to $\left(
Y,Z,U,K\right)  \in S_{\mathbb{G}}^{2}\times H_{\mathbb{G}}^{2}\times
L^{2}(\lambda)\times\mathcal{K}^{2}$, such that%
\[
\left\{
\begin{array}
[c]{ll}%
Y_{t} & =\xi^{2}+%
{\textstyle\int_{t}^{T}}
g^{2}(s,Y_{s},Z_{s},\mathbb{E}[Y_{s+\delta}|\mathcal{G}_{s}],U_{s})ds-%
{\textstyle\int_{t}^{T}}
Z_{s}dW_{s}\\
& \qquad-%
{\textstyle\int_{t}^{T}}
U_{s}dM_{s}+%
{\textstyle\int_{t}^{T}}
dK_{s},\text{\ \qquad}t\in\left[  0,T\right]  ,\\
Y_{t} & =\xi^{2},\text{\qquad}t\geq T,\\
Z_{t} & =U_{t}=0,\text{\qquad}t>T.
\end{array}
\right.
\]
By the uniqueness of the solution (Theorem \ref{ex-uni}), we have that
$Y_{t}=Y_{t}^{2},$ a.s.

\noindent Since for all $t,$ $Y_{t}^{1}\geq Y_{t}^{3}\geq Y_{t}^{4}\geq...\geq
Y_{t},$ a.s. it hold immediately for a.a. $t$%
\[
Y_{t}^{1}\geq Y_{t}^{2},\text{ a.s.}%
\]
$\qquad\qquad\qquad\qquad\qquad\qquad\qquad\qquad\qquad\qquad\qquad
\qquad\qquad\square$ \newline

\section{RABSDE with Default and Optimal Stopping}

We recall here a connection between RABSDE and optimal stopping problems. The
following result is essentially due to El Karoui \textit{et al} \cite{elkarui}
under the Brownian filtration and to {Ø}ksendal and Zhang \cite{oz}:

\begin{definition}
$\bullet$ Let $F:\Omega\times\lbrack0,T]\times\mathbb{R}^{5}\rightarrow
\mathbb{R}$, be a given function such that:\newline

$\bullet$ $F$ is $\mathbb{G}$-adapted and $\mathbb{E[}%
{\textstyle\int_{0}^{T}}
|F(t,0,0,0,0,0)|^{2}dt]<0.$

$\bullet$ Let $S_{t}$ be a given $\mathbb{G}$-adapted continuous process such
that $\mathbb{E[}\underset{t\in\lbrack0,T]}{\sup}S_{t}^{2}]<\infty$.

$\bullet$ The terminal value $\xi\in L^{2}\left(  \Omega,\mathcal{G}%
_{T}\right)  $ is such that $\xi\geq S_{T}$ a.s.

We say that a $\mathbb{G}$- adapted triplet $\left(  Y,Z,K\right)  $ is a
solution of the reflected ABSDE with driver $F$, terminal value $\xi$ and the
reflecting barrier $S_{t}$ under the filtration $\mathbb{G}$, if the following hold:

\begin{enumerate}
\item $\mathbb{E[}%
{\textstyle\int_{0}^{T}}
|F(s,Y_{s},Z_{s},\mathbb{E[}Y_{s+\delta}|\mathcal{G}_{s}],\mathbb{E[}%
Z_{s+\delta}|\mathcal{G}_{s}],U_{s})|^{2}dt]<\infty$,

\item $Y_{t}=\xi+%
{\textstyle\int_{t}^{T}}
F(s,Y_{s},Z_{s},\mathbb{E[}Y_{s+\delta}|\mathcal{G}_{s}],\mathbb{E[}%
Z_{s+\delta}|\mathcal{G}_{s}],U_{s})ds-%
{\textstyle\int_{t}^{T}}
dK_{s}-%
{\textstyle\int_{t}^{T}}
Z_{s}dW_{s}-%
{\textstyle\int_{t}^{T}}
U_{s}dM_{s},t\in\left[  0,T\right]  ,$\newline or, equivalently,\newline%
$Y_{t}=\mathbb{E}[\xi+%
{\textstyle\int_{t}^{T}}
F(s,Y_{s},Z_{s},\mathbb{E[}Y_{s+\delta}|\mathcal{G}_{s}],\mathbb{E[}%
Z_{s+\delta}|\mathcal{G}_{s}],U_{s})ds-%
{\textstyle\int_{t}^{T}}
dK_{s}|\mathcal{G}_{t}],t\in\left[  0,T\right]  ,$

\item $K_{t}$ is nondecreasing, $\mathbb{G}$-adapted, càdlàg process with $%
{\textstyle\int_{0}^{T}}
(Y_{t}-S_{t})dK_{t}^{c}=0$ and $\triangle K_{t}^{d}=-\triangle Y_{t}%
\mathbf{1}_{\{Y_{t^{-}}=S_{t^{-}}\}}$, where denote the continuous and
discontinuous parts of $K$ respectively$,$

\item $Y_{t}\geq S_{t}$ a.s., $t\in\lbrack0,T].$
\end{enumerate}

\vskip0.3cm
\end{definition}

\begin{theorem}
For $t\in\lbrack0,T]$ let $\mathcal{T}_{[t,T]}$ denote the set of all
$\mathbb{G}$-stopping times $\tau:\Omega\mapsto\lbrack t,T].$\newline Suppose
$\left(  Y,Z,U,K\right)  $ is a solution of the RABSDE above.
\end{theorem}

\begin{description}
\item[(i)] Then $Y_{t}$ is the solution of the optimal stopping problem%
\[%
\begin{array}
[c]{c}%
Y_{t}=\underset{\tau\in\mathcal{T}_{[t,T]}}{ess\sup}\quad\{\mathbb{E}[%
{\textstyle\int_{t}^{\tau}}
F(s,Y_{s},Z_{s},\mathbb{E[}Y_{s+\delta}|\mathcal{G}_{s}],\mathbb{E[}%
Z_{s+\delta}|\mathcal{G}_{s}],U_{s})ds\\
+S_{\tau}\mathbf{1}_{\tau<T}+\xi\mathbf{1}_{\tau=T}|\mathcal{G}_{t}]\},\quad
t\in\left[  0,T\right]  .
\end{array}
\]

\item[(ii)] Moreover the solution process $K(t)$ is given by%
\[%
\begin{array}
[c]{c}%
K_{T}-K_{T-t}=\underset{s\leq t}{\max}\{\xi+%
{\textstyle\int_{T-s}^{T}}
F(r,Y_{r},Z_{r},\mathbb{E[}Y_{r+\delta}|\mathcal{G}_{r}],\mathbb{E[}%
Z_{r+\delta}|\mathcal{G}_{r}],U_{r})dr\\
-%
{\textstyle\int_{T-s}^{T}}
Z_{r}dB_{r}-S_{T-s}\},\quad t\in\left[  0,T\right]  ,
\end{array}
\]
where $x^{-}=\max(-x,0)$ and an optimal stopping time $\hat{\tau}_{t}$ is
given by%
\begin{align*}
\hat{\tau}_{t}:  &  =\inf\{s\in\lbrack t,T],Y_{s}\leq S_{s}\}\wedge T\\
&  =\inf\{s\in\lbrack t,T],K_{s}>K_{t}\}\wedge T.
\end{align*}

\item[(iii)] In particular, if we choose $t=0$ we get that
\begin{align*}
\hat{\tau}_{0}: &  =\inf\{s\in\lbrack0,T],Y_{s}\leq S_{s}\}\wedge T\\
&  =\inf\{s\in\lbrack0,T],K_{s}>0\}\wedge T
\end{align*}
solves the optimal stopping problem
\[%
\begin{array}
[c]{c}%
Y_{0}=\sup_{\tau\in\mathcal{T}_{[0,T]}}\mathbb{E}[{%
{\textstyle\int_{0}^{\tau}}
}F(s,Y_{s},Z_{s},\mathbb{E[}Y_{s+\delta}|\mathcal{G}_{s}],\mathbb{E[}%
Z_{s+\delta}|\mathcal{G}_{s}],U_{s})ds\\
+S_{\tau}\mathbf{1}_{\tau<T}+\xi\mathbf{1}_{\tau=T}],\quad t\in\left[
0,T\right]  .
\end{array}
\]

\end{description}

\textbf{Acknowledgement}. We would like to thank Professors Bernt Øksendal and
Shiqi Song for helpful discussions.

\end{document}